\documentclass[12pt]{amsart}
\usepackage{latexsym, ifthen,amsfonts,amsmath,amssymb}
\usepackage[dvips]{graphics}
\usepackage{epsfig,color}
 \normalbaselines
 \newcommand{\newsection}[1]
 {\subsection{#1}\setcounter{theorem}{0} \setcounter{equation}{0}
 \par\noindent}
 \renewcommand{\theequation}{\arabic{subsection}.\arabic{equation}}
 \renewcommand{\thesubsection}{\arabic{subsection}}
 \newtheorem{theorem}{Theorem}
 \renewcommand{\thetheorem}{\arabic{subsection}.\arabic{theorem}}
 \newtheorem{lemma}[theorem]{Lemma}
 \newtheorem{corr}[theorem]{Corollary}
 
 \newtheorem{proposition}[theorem]{Proposition}
 \newtheorem{deff}[theorem]{Definition}
 
 \newcommand{\bth}{\begin{theorem}}
 \newcommand{\ble}{\begin{lemma}}
 \newcommand{\bcor}{\begin{corr}}
 \newcommand{\bdeff}{\begin{deff}}
 \newcommand{\bprop}{\begin{proposition}}
 \renewcommand{\eth}{\end{theorem}}
 \newcommand{\ele}{\end{lemma}}
 \newcommand{\ecor}{\end{corr}}
 \newcommand{\edeff}{\end{deff}}
 \newcommand{\beq}{\begin{equation}}
 \newcommand{\eeq}{\end{equation}}
 \newcommand{\eprop}{\end{proposition}}

  \newcommand{\B}{{\emph{B}}}

 \renewcommand{\Pi}{\varPi}

 \renewcommand{\epsilon}{\varepsilon}

 \newcommand{\R}{{\mathbb{R}}}

  \newcommand{\rng}{\rangle}
   \newcommand{\lng}{\langle}

\newcommand{\BB}{\mathcal{B}}

\renewcommand{\S}{\Sigma}
\renewcommand{\SS}{\Sigma}
\newcommand{\DD}{\mathcal{D}}
\newcommand{\SP}{\text{ }}

\begin{document}
\title {Multi-valued graphs in embedded constant mean curvature
disks.} \maketitle

 \noindent GIUSEPPE TINAGLIA, Department of
Mathematics,
Johns Hopkins University, 3400 North Charles Street, 404 Krieger Hall, Baltimore, MD 21218-2686.\\
e-mail: tinaglia@math.jhu.edu
\begin{abstract}
In this paper we prove that an embedded constant mean curvature
disk with curvature large at a point contains a multi-valued graph
around that point on the scale of $|A|^2$. This generalizes
Colding and Minicozzi's result for minimal surfaces.
\end{abstract}

\section*{Introduction}
\begin{figure}[htbp]
    \setlength{\captionindent}{4pt}
    \begin{minipage}[t]{0.5\textwidth}
    \centering\input{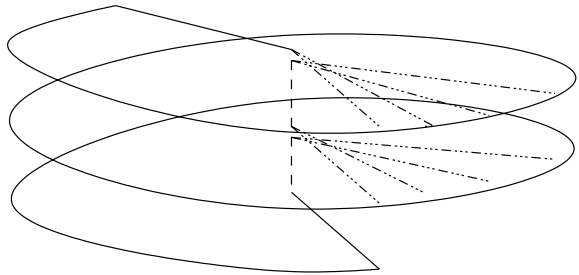}
   \caption{Half of the the Helicoids}\label{f:f1}
    \end{minipage}
\end{figure}

In this paper we prove that an embedded and simply connected
constant mean curvature (CMC) surface with curvature large at a
point contains a multi-valued graph around that point on the scale
of $|A|^2$, where $|A|^2$ is the norm squared of the second
fundamental form.
%
More precisely, our main result is the following:

\bth \label{CoMi} Given $N \in \mathbb {Z}_+,\SP\omega>1$ and
$\epsilon>0$, there exist $C=C(N,\omega,\epsilon)>0$, $H>0$ and $\bar{l}>1$ so:\\
Let $\Sigma\subset\R^3$ be an embedded and simply connected
constant mean curvature equal to $h$ surface. If
$|h|<\frac{H}{r_0}$ and
$$\sup_{\Sigma\cap B_{r_0\bar{l}}(0)}|A|^2\leq 4C^2r_0^{-2}=4|A|^2(0)$$
 for some $r_0>0$, then $\Sigma$ (after a rotation) contains an $N$-valued
graph over $D_{\omega \bar{R}}\backslash D_{\bar{R}}$ where
$\bar{R}<\frac{r_0}{\omega}$ (with gradient $\leq \epsilon$ and
$\text{\textsl{dist}}_{\Sigma}(0,\Sigma_g)\leq4\bar{R}$). \eth
Roughly speaking, to contain a multi-valued graph
(Def.~\ref{multigraph} in this paper) means that locally the
surface spirals like a helicoid, Fig. \ref{f:f1}.
The helicoid is a minimal surface parameterized in the following
way
$$(s\sin t,s\cos t,t) \quad\quad \text{where } (s,t) \in \R^2 .$$

Our result is a generalization of Colding and Minicozzi's result
\cite[\textbf{Theorem 0.4.}]{CM2} (Thm. \ref{result} in this
paper) which is a key ingredient in their series of papers
\cite{CM1,CM2,CM3,CM4} that dealt with the structure of embedded
minimal disks. We prove that under equivalent local conditions an
embedded CMC disk contains a multi-valued graph as well. For a
minimal surface, Colding and Minicozzi were able to extend the
multi-valued graph that forms locally, all the way up to the
boundary \cite{CM1}. It is not known if the same can be done for
CMC surfaces.

The proof is by contradiction using a compactness argument. The
idea is the following: Assuming that Theorem \ref{CoMi} is false,
we build a sequence $\S_n$ of embedded CMC disks where each disk
satisfies the hypotheses of the theorem with $C$ fixed large and
$H=\frac{1}{n}$ but does not contain a $N$-valued graph.
We prove that $\S_n$ converges $C^2$ to a minimal surface
$\Sigma_\infty$ which contains an $N$-valued graph.
\begin{deff}\label{convergence}
A sequence $\S_n$ of surfaces converges to a surface $\S_\infty$
in the $C^k$ topology if at any point $p\in\S_\infty$ each $\S_n$
is locally (near $p$) a graph over the tangent space
$T_p\S_\infty$ and the graph of $\S_n$ converges to the graph of
$\S_\infty$ in the usual $C^k$ topology.
\end{deff}
We will also consider sequences that converge with multiplicity.
This means that we allow more than one graph in the previous
definition.

Essentially we show that $\S_n$ comes as close as we want to its
limit and that the limit is an embedded minimal disk which
contains an $N$-valued graph because of \textbf{Theorem 0.4.} in
\cite{CM2}, therefore so do the CMC surfaces. To create the
$N$-valued graph in the CMC sequence we basically push the
multi-valued graph from the minimal surface onto $\S_n$. This
contradiction proves the theorem.

The difficult part of the proof is to show that the limit is both
an embedded surface and simply connected and not, for instance, a
minimal lamination or a minimal surface which is not simply
connected. Some sort of $C^2$ convergence follows in a standard
way from the bound on the curvature and trivially, since we are
assuming that the mean curvature goes to zero, the limit is a
minimal object. To assure that the limit is embedded and simply
connected we need a uniform upper bound on the number of graphs
over the tangent plane $T_p\S_\infty$. In order to obtain this
uniform upper bound, we investigate the strong stability for a
constant mean curvature surface to find out when a CMC surface,
which is already a critical point for a certain area functional,
is an actual minimum.

First, we prove that, under certain conditions, if two CMC
surfaces are close and disjoint, they are almost-stable ("almost a
minimum").

Second, we rule out the possibility that $\S_n$ contains a large,
almost-stable domain, for $n$ large.

Third, we show that if there is not a uniform upper bound on the
number of pieces, then two large pieces of $\S_n$ are close and
disjoint, creating a large almost-stable domain and giving the
contradiction.

Once the uniform upper bound on the number of pieces is obtained,
the convergence to an embedded minimal surface follows. We have to
use some topological results to prove that the limit minimal
surface is simply connected.

In the first section we provide a short overview of constant mean
curvature surfaces. In the second section we describe what a
multi-valued graph is and go over the hypotheses of the main
result. We also take a closer look at the proof. In the third
section we deal with the $\delta$-stability for CMC surfaces, and
give a criteria to find $\delta$-stable domains in CMC surfaces.
In the fourth section we show how, because of the upper-bound on
$|A|^2$, our CMC disk is "uniformly locally flat" and we give a
criteria to find large $\delta$-stable domains. In the fifth
section we prove that a large $\delta$-stable domain cannot be
contained in $\S_n$ for $n$ large and how this gives an upper
bound on the number of graphs. In the sixth section we prove that
the limit is an embedded minimal disk and from that we build a
multi-valued graph in $\S_n$.

We actually prove the result when $r_0$ in Theorem \ref{CoMi} is
fixed and equal to one. The main result will follow by rescaling
and in Appendix B we describe the rescaling argument for constant
mean curvature graphs. In Appendix A we provide examples of CMC
surfaces containing arbitrary large multi-valued graphs.

\newsection{Constant Mean Curvature Surfaces}\label{one}

This section is a short review of general properties of CMC
surfaces.

Let $\Sigma\subset\R^3$ be a 2-dimensional smooth orientable
surface (possibly with boundary) with unit normal $N_{\Sigma}$.
Given a function $\phi$ in the space $C^{\infty}_0(\Sigma)$ of
infinitely differentiable (i.e., smooth), compactly supported
functions on $\Sigma$, consider the one-parameter variation
$$ \Sigma_{t,\phi}=\{x+t\phi(x)N_{\Sigma}(x)|x\in\Sigma\} $$
and let $A(t)$ be the area functional,
$$ A(t)=Area(\Sigma_{t,\phi}).$$
The so-called first variation formula of area is the equation
(integration is with respect to $d$area)
\beq A'(0)=\int_{\Sigma}\phi H,\eeq
where $H$ is the mean curvature of $\Sigma$. When $H$ is constant
the surface is said to be a \textit{constant mean curvature} (CMC)
surface \cite{Ken} and it is a critical point for the area
functional restricted to those variations which preserve the
\textit{enclosed volume}, in other words $\phi$ must satisfy the
condition,
$$\int_{\Sigma}\phi=0.$$
In general, if $\Sigma$ is given as graph of a function $u$ then
\beq H=\text{div} \left (\frac{\nabla u}{\sqrt{1+|\nabla u|^2}}
\right ) .\eeq
Therefore, when $H$ is constant $u$ satisfies a quasi-linear
differential equation. In the particular case where the mean
curvature $H$ is identically zero the surface $\Sigma$ is said to
be a \textit{minimal} surface \cite{Oss, CM}. Concrete examples of
constant mean curvature surfaces are spheres, cylinders and
Delauney surfaces.

Let $u_1$, $u_2$ be CMC graphs over $D_r(0)$ and assume that they
have the same constant mean curvature ($H_{u_1}=H_{u_2}$), the
same orientation ($\lng N_1,N_2 \rng >0$), and that $u_1-u_2>0$
then \cite[Lemma 1.17]{CM} $u_1-u_2=v$ is a positive solution of
\begin{equation}
\text{div}A_{ij}\nabla v+b\nabla v=0 
\end{equation}
where $A_{ij}$ and $b$ depend on $\nabla u_1$ and $\nabla u_2$.
Moreover, if $|\nabla u_1|$ and $|\nabla u_2|$ are sufficientely
small, we have the Harnack type inequality
\begin{equation}\label{CO}
\sup_{B_\frac{r}{2}(0)}(u_1-u_2)\le C_0(u_1(0)-u_2(0)).
\end{equation}
Notice from Fig. \ref{circles} that the condition $\lng
N_1,N_2\rng>0$ on the orientation is necessary. As it is shown in
Fig. \ref{circles} the two spherical caps have the same constant
mean curvature since they have the same radius. However, even if
$u_1(0)-u_2(0)=0$, it is clear that $\sup (u_1-u_2)>0$ in any
neighborhood of 0 and therefore that \eqref{CO} does not follow.
\begin{figure}[htbp]
\begin{center}
    \centering\input{fig1a.pstex_t}
    \end{center}
    \caption{}
   \label{circles}
\end{figure}
%
In general, let $k_1, k_2$ be the principal curvatures on
$\Sigma$, then $H=k_1+k_2$; $|A|^2=k_1^2+k_2^2$ is the norm
squared of the second fundamental form. Since the Gaussian
curvature $K_\Sigma$ is equal to the product of the principal
curvatures $ k_1k_2$, we have the Gauss equation, that is
\beq \label{gauss} H^2=|A|^2+2K_\Sigma .\eeq
From \eqref{gauss} it is clear why when $H$ is constant, in
particular when it is small and even better when it is zero,
talking about the Gaussian curvature or talking about the norm of
the second fundamental form squared is almost equivalent.

\newsection{Multi-valued graphs in CMC surfaces}\label{proof}

In this section we discuss the result and explain the necessity of
the hypotheses. We also take a closer look at how the proof goes.

\begin{figure}[htbp]
    \setlength{\captionindent}{4pt}
    \begin{minipage}[t]{0.5\textwidth}
    \centering\input{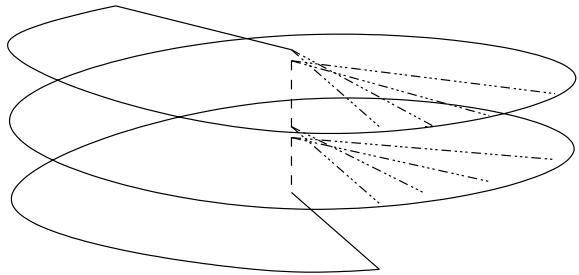}
   \caption{Half of the the Helicoids}\label{fig1a}
    \end{minipage}
\end{figure}

This is what Colding and Minicozzi proved:
\bth \label{result}\cite[\textbf{Theorem 0.4.}]{CM2} Given $N \in
\mathbb {Z}_+,\SP\omega>1$ and $\epsilon>0$, there exist
$C=C(N,\omega,\epsilon)>0$ so:\\
Let $0\in\Sigma\subset B_{R}\subset\R^3$ be an embedded minimal
disk such that $\partial \Sigma \subset B_{R}$. If
$$\sup_{\Sigma\cap B_{r_0}}|A|^2\leq 4C^2r_0^{-2}\text{ and
}|A|^2(0)=C^2r_0^{-2}$$ for some $0<r_0<R$, then there exists
$\bar{R}<\frac{r_0}{\omega}$ and (after a rotation) an $N$-valued
graph $\Sigma_g\subset\Sigma$ over $D_{\omega \bar{R}}\backslash
D_{\bar{R}}$ with gradient $\leq \epsilon$ and
$\text{\textsl{dist}}_{\Sigma}(0,\Sigma_g)\leq4\bar{R}$. \eth
\begin{deff}[Multigraph]\label{multigraph}
Let $D_r$ be the disk in the plane centered at the origin and of
radius $r$ and let $\mathcal{P}$ be the universal cover of the
punctured plane $\mathbb{C}\backslash{0}$ with global coordinates
$(\rho , \theta)$ so $\rho >0$ and $\theta\in\mathbb{R}$. An
$N$-valued graph of a function $u$ on the annulus $ D_s \backslash
D_r$ is a single valued graph over $\{(\rho ,\theta )| r\leq \rho
\leq s,|\theta |\leq N\pi \}$.
\end{deff}
When dealing with multi-valued graphs, the surface to keep in mind
is the helicoid, Fig. \ref{fig1a}. A parametrization of the
helicoid that illustrates the existence of such an $N$-valued
graph is the following
$$(s\sin t,s\cos t,t) \quad\quad \text{where } (s,t) \in \R^2 .$$
It is easy to see that it contains the $N$-valued graph $\phi$
defined by
$$\phi(\rho,\theta)=\theta \quad\quad \text{where }(\rho ,\theta)\in \R^+ \backslash 0 \times [-N\pi,N\pi].$$
In fact the helicoid is a minimal surface. In Appendix A we
provide examples of CMC surfaces containing
arbitrary large multi-valued graphs.\\

What we are about to prove is not exactly Theorem \ref{CoMi}. We
prove the result when $r_0$ in Theorem \ref{CoMi} is equal to one
and hence the curvature is bounded in a ball of radius $\bar{l}$.
We will discuss and determine $\bar{l}$ in Section 5. Theorem
\ref{CoMi} will follow by rescaling and we will describe the
rescaling argument in Appendix B.

This is the new statement:
\bth \label{main2} For each $N \in \mathbb {Z}_+$, $\omega>1$ and
$\epsilon>0$ there exist $H>0$, $C(N,\omega,\epsilon)>0$ and
$\bar{l}>1$ so:\\
Let $0\in\Sigma\subset B_{\bar{l}}(0)\subset\R^3$ be an embedded
and simply connected constant mean curvature surface equal to $h$
(embedded CMC disk) such that $|h|\leq H$ and $\partial \Sigma
\subset
\partial B_{\bar{l}}(0)$. If
$$\sup_{\Sigma\cap B_{\bar{l}}(0)}|A|^2\leq 4C^2=4|A|^2(0)$$ then there exists $\bar{R}<\frac{1}{\omega}$ and
(after a rotation) an $N$-valued graph $\Sigma_g\subset\Sigma$
over $D_{\omega \bar{R}}\backslash D_{\bar{R}}$ (with gradient
$\leq \epsilon$ and
$\textsl{dist}_{\Sigma}(0,\Sigma_g)\leq4\bar{R}$). \eth

The constant $C(N,\omega,\epsilon)$ is essentially the same
constant that Colding and Minicozzi used.

We can only prove that a multi-valued graph exists substantially
far away from the boundary, that is in a ball of radius one while
the boundary of the surface is contained in the boundary of a ball
of radius $\bar{l}>1$. For a minimal surface, Colding and
Minicozzi were able to extend the multi-valued graph that forms
locally, all the way up to the boundary \cite{CM1}. It is not
known if the same can be done for CMC surfaces.

Thanks to the upper bound on the second fundamental form, the
surface is "uniformly locally flat" and the $C^2$ convergence
follows. Moreover, $\sup_{\Sigma\cap B_{\bar{l}}(0)}|A|^2\leq
4C^2$ together with the Gauss equation \eqref{gauss} gives a lower
bound for the Gaussian curvature,
\begin{equation}\label{bishop1}
K_\Sigma\geq -4C^2=2G.
\end{equation}
This lower bound implies an upper bound on the area of the
intrinsic balls, Theorem \ref{bishop}.

%
%

%
What follows is a short sketch of the proof. The proof is by
contradiction. Assuming that Theorem \ref{main2} is false we have
the following:
\begin{quote}
Given $C(N,\omega,\epsilon)$ as in Theorem \ref{result}, for any
$h>0$ there exists an embedded and simply connected constant mean
curvature equal $h$ surface $\Sigma_h$ that does not contain an
$N$-valued graph $\Sigma_g\subset\Sigma$ over $D_{\omega
\bar{R}}\backslash D_{\bar{R}}$ for any $\bar{R}<\frac{1}{\omega}$
but such that
$$0\in\Sigma\subset B_{\bar{l}}(0)\subset\R^3\text{, }\partial \Sigma
\subset
\partial B_{\bar{l}}(0)\text{ and }\sup_{\Sigma\cap B_{\bar{l}}}|A|^2\leq 4C^2=4|A|^2(0).$$
\end{quote}
We want to show that this cannot be true. Let us take a sequence
of $\S_n$ as above with $h=\frac{1}{n}$. The constant mean
curvature of $\S_n$ goes to zero but none of the elements in the
sequence contain an $N$-valued graph. Fixed $\bar{\epsilon}>0$, we
consider a new sequence $\S_n'$ where $\S_n'$ is the connected
component of $\S_n\cap B_{\bar{l}-\bar{\epsilon}}(0)$ that
contains 0. Given that $|A|^2$ is bounded and we are slightly away
from the boundary
there exists $r>0$ so:\\
$\S_n'$ can be covered by a finite number of balls, $B_r(x^n_i)$
where $x^n_i\in\S_n'$, such that in each ball $\Sigma_n'\cap
B_r(x^n_i)$ looks like graphs $u^j_n$ over the tangent plane
$T_{x^n_i}\S_n$. The radius $r$ and the number of balls will be
independent of $n$. Going to a subsequence, we can assume that
$x^n_i$ converges to a certain $x_i$ and that $T_{x^n_i}\S_n$
converges to a certain $T_{x_i}\S_{\infty}$. At this point we are
able to extract, by using Arzela-Ascoli, a subsequence $u^j_n$
that converges uniformly to a graph $u^j_\infty$. These CMC graphs
satisfy the following partial differential equation
$$\frac{1}{n}=\text{div}\left (\frac{\nabla u^j_n}{\sqrt{1+|\nabla u^j_n|^2}}\right ). $$
Therefore, using Schauder theory \cite{GiTru} and the fact that
$\frac{1}{n}$ goes to zero, we can prove that $u^j_n$ converges
$C^2$ to $u^j_\infty$ and that the latter is a minimal graph.

\begin{figure}[htbp]
    \setlength{\captionindent}{4pt}
    \begin{minipage}[t]{0.5\textwidth}
    \centering\input{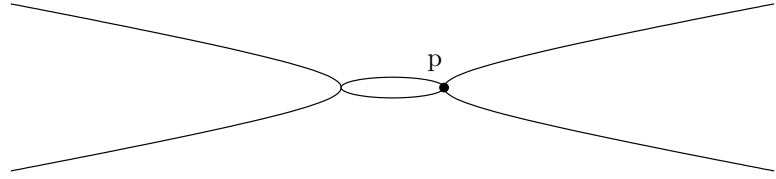}
   \caption{Catenoid}\label{catenoid}
    \end{minipage}
\end{figure}

Unfortunately, we need more to prove the global properties
required. The limit object contains a multi-valued graph if it is
an embedded and simply connected minimal surface. We have not
ruled out the possibility that the number $j$ of graphs $u^j_n$
goes to infinity as $n$ goes to infinity and in the limit that
could give an infinite number of minimal graphs. As a consequence
the limit would not necessarily be a surface but it could be a
lamination. Another possibility is that the limit is not simply
connected, for instance it could be a catenoid, Fig.
\ref{catenoid}.
Rescaling the catenoid in Fig. \ref{catenoid} the curvature at $p$
becomes very large and yet the catenoid would not contain a
multi-valued graph.

What we show is that the number of graphs is uniformly bounded if
we stay substantially away from the boundary. This is because to
prove this uniform upper-bound on the number of graphs we have to
work with large geodesic balls and to assure that they exist, we
need to move substantially away from the boundary. We need to be
working in the unit ball and keep the boundary of the surface on a
substantially bigger ball. More precisely, we build another
subsequence $\S^1_n$ where $\S^1_n$ is the connected component of
$\S_n\cap B_1(0)$ that contains 0, Fig. \ref{three}.
\begin{figure}[htbp]
    \setlength{\captionindent}{4pt}
    \begin{minipage}[t]{0.5\textwidth}
    \centering\input{fig3a.pstex_t}
   \caption{}
   \label{three}
    \end{minipage}
\end{figure}
$\S^1_n$ is also simply connected. If it was not simply connected
there would exist $B_r(0)$, $1<r<\bar{l}$ such that $\S_n$ is
tangent to $\partial B_r(0)$ and locally inside $B_r(0)$. This is
a contradiction for $H(n)<\frac{1}{2\bar{l}}$, Fig. \ref{sc}.
\begin{figure}[htbp]
    \setlength{\captionindent}{4pt}
    \begin{minipage}[t]{0.5\textwidth}
    \centering\input{fig2d.pstex_t}
  \caption{}
   \label{sc}
    \end{minipage}
\end{figure}
If we restrict our attention to $\S^1_n$ we have a uniform upper
bound on the number of graphs and it follows that $\S^1_n$, not
the whole $\S_n$, converges to an embedded minimal disk. Once we
have that $\S_n^1$ converges to an embedded minimal disk, we prove
that $\S_n^1$, and therefore $\S_n$, contains a multi-valued
graph.

From now on, even if the results can often be stated more
generally, $\S$ will be a CMC surface satisfying the hypotheses of
Theorem $\ref{main2}$. $\S'$ will be the connected component of
$\S\cap B_{\bar{l}-\bar{\epsilon}}(0)$ containing 0. $\S^1$ will
be the connected component of $\S\cap B_1(0)$ containing 0. We
will also assume the mean curvature to be as small as we need, in
particulary bounded.

\newsection{$\delta$-stability}

This section consists of standard results about CMC surfaces and
stability.

Let $A$ be the area functional described in Section \ref{one}; we
showed that $A'(0)=\int_{\Sigma}\phi H$. A computation shows that
if $\Sigma$
is a CMC surface then%
\beq A''(0)=-\int_{\Sigma}\phi L_{\Sigma}\phi, \quad\quad
\text{where } L_{\Sigma}\phi=\Delta_{\Sigma}\phi +|A|^2\phi \eeq%
is the second variational operator. Here $\Delta_{\Sigma}$ is the
intrinsic Laplacian on $\Sigma$. A CMC surface $\Sigma$ is said to
be (strongly) stable if
\beq \label{stable} A''(0)\geq 0 \quad \quad \text{for all } \phi
\in C_0^{\infty}(\Sigma).\eeq %
Applying Stokes' theorem to \eqref{stable} shows that $\Sigma$ is
stable if and only if
$$ \int_{\Sigma}|A|^2\phi^2\leq
\int_{\Sigma}|\nabla\phi|^2, \quad\quad \text{for all } \phi \in
C_0^{\infty}(\Sigma) $$
and that allows us to define $\delta$-stability, namely $\Sigma$
is said to be $\delta$-stable if
\beq (1-\delta)\int_{\Sigma}|A|^2\phi^2\leq
\int_{\Sigma}|\nabla\phi|^2, \quad\quad \text{for all } \phi \in
C_0^{\infty}(\Sigma).\eeq

In the following lemma we establish a relation between a CMC
surface and a CMC normal variation of it that does not change the
mean curvature.

\ble\label{stability} There exists $\delta_1>0$ so: If
$\delta<\delta_1$, $\Sigma$ is a CMC surface and $u$ is a positive
solution of the CMC graph equation over $\Sigma$ (i.e.
$\Sigma^u:=\{x+u(x)N_{\Sigma}(x)|x\in\Sigma\}$ is CMC) such that
$|H_{\Sigma^u}|=|H_\Sigma|$, $\lng N_{\Sigma^u},N_\Sigma \rng
\geq0$, $|u||A|$ and $|\nabla{u}|\leq\delta$ then
$\triangle{u}+u|A|^2=o(\delta^2)$. \ele
\begin{proof}
In general
$$H_{\Sigma^u}=H_\Sigma+\frac{1}{2}(\triangle{u}+u|A|^2)+o(|u|^2,|\nabla u|^2).$$
The condition $\lng N_{\Sigma^u},N_\Sigma \rng \geq0$ is a
condition on the orientation that implies $H_{\Sigma^u}=H_\Sigma$
and the lemma follows.
\end{proof}
The existence of a positive solution of $Lu=0$ where $L$ is
$\Delta +|A|^2$ would imply $A''(0)\geq 0$ for all $\phi \in
C_0^{\infty}(\Sigma)$. In the following lemma we show that if
there exists a positive function $u$ which is "almost" a solution,
then $A''(0)$ is "almost" non-negative for all $\phi \in
C_0^{\infty}(\Sigma)$, that is, almost-stable.
\ble\label{stability2} Let $\Omega$ be a domain and $u$ be a
positive function in $C^2(\Omega)$ such that %
\beq \label{omega} \Delta u \leq -(1-\delta)|A|^2u\eeq%
then $\Omega$ is  $ \delta$-stable.\ele
\begin{proof}
Set $w=\log{u}$ and let $\Phi$ be any
compactly supported function on $\Omega$. We have %
$$\text{div}(\nabla{w})=\text{div}(\frac{\nabla{u}}{u})=\frac{\triangle{u}}{u}-\frac{|\nabla{u}|^2}{u^2}= \frac{\triangle{u}}{u}-|\nabla{w}|^2.$$\\
Applying Stokes theorem to $\text{div}(\Phi^2\nabla{w})$ gives
$$0=\int \text{div}(\Phi^2\nabla{w})=\int \Phi^2\Delta w + \int\lng\nabla{\Phi}^2,\nabla{w}\rng.$$ Using Cauchy-Shwarz and the absorbing inequality gives
$$\int \lng\nabla{\Phi}^2,\nabla{w}\rng\leq\int|\nabla{\Phi}|^2 + \int
\Phi^2|\nabla{w}|^2.$$
Eventually,
$$\int(-\frac{\triangle{u}}{u}+|\nabla{w}|^2) \Phi^2 \leq \int|\nabla{\Phi}|^2 + \int
\Phi^2|\nabla{w}|^2.$$ Applying \eqref{omega} we get
$$(1-\delta)\int_{\Sigma}|A|^2\phi^2\leq
\int_{\Sigma}|\nabla\phi|^2 .$$
\end{proof}
%

Lemma \ref{stability} and Lemma \ref{stability2} give a first
criteria to find almost stable domains in a constant mean
curvature surface.
\begin{corr}\label{criteria1}
There exists $\delta_3>0$ so: If $\delta<\delta_3$, $\Sigma$ is a
CMC surface and $u$ is a positive solution of the CMC graph
equation over $\Sigma$ such that $|H_{\Sigma^u}|=|H_\Sigma|$,
$\lng N_{\Sigma^u},N_\Sigma\rng\geq0$, $|u||A|$ and
$|\nabla{u}|\leq\delta$ then $\S$ is $\delta$-stable.
\end{corr}
\newsection{The upper bound on $|A|^2$}

In this section we use the upper bound on $|A|^2$ to generalize
some standard local results regarding CMC surfaces. We prove a
criteria to find large pieces of $\S'$ which are graph over other
pieces, creating large almost stable CMC domains.

Let us define
\beq \SS_{x,R} \text{ as the component of } B_R(x)\cap\SS \text{
that contains } x, \eeq
\beq \BB_R(x):=\{y\in \S \text{ such that }
\text{dist}_\S(y,x)<R\}\eeq i.e., the geodesic ball of radius $R$
centered at $x$,
\beq \DD_r(x):=\{x'\in T_x\S \text{ such that } |x-x'|<r\}.\eeq

In what follows we are about to explain why in a CMC surface with
bounded $|A|^2$ everything looks graphical--what we have been
calling "uniformly locally flat." Integrating
$|\nabla\text{dist}_{S^2}(\textbf{n}(x),\textbf{n})|\leq|A|$ on
geodesics gives
\begin{equation}\label{eq1}
\sup_{x'\in\BB_s(x)}\text{dist}_{S^2}(\textbf{n}(x),\textbf{n})\leq
s\sup_{\BB_s(x)}|A|.
\end{equation}
By \eqref{eq1}, we can choose $0<\rho<\frac{1}{4}$ so: If
$\BB_{2s}(x)\subset\SS$, $s\sup_{\BB_s(x)}|A|\leq4\rho_2$, and
$t\leq s$ then the component $\SS_{x,t}$ of $B_t(x)\cap\SS$ with
$x\in\SS_{x,t}$ is a graph over $T_x\SS$ with gradient
$\leq\frac{t}{s}$ and
\begin{equation}\label{eq2}
1\geq\inf_{x'\in\BB_2s(x)}\frac{|x'-x|}{\text{dist}_{\SS}(x,x')}>\frac{9}{10}.
\end{equation}
One consequence is that if $t\leq s$ and we translate $T_x\SS$ so
that $x\in T_x\SS$, then
\begin{equation}\label{eq3}
\sup_{x'\in\BB_t(x)}|x'-T_x\SS|\leq \frac{t^2}{s}.
\end{equation}
\begin{figure}[h]
    \setlength{\captionindent}{4pt}
    \begin{minipage}[t]{0.5\textwidth}
    \centering\input{fig5.pstex_t}
  \caption{}
   \label{graph}
    \end{minipage}
\end{figure}%
As a consequence of \eqref{eq1}, \eqref{eq2}, \eqref{eq3} and the
fact that in this paper we are assuming $\sup_{\SS}|A|<C$ we can
clearly choose $0<\bar{\rho}<4\frac{\rho}{C}$ so: Given
$t<\bar{\rho}$ and $x\in\SS$ then $$\S_{x,t}\text{ is a graph over
} T_x\SS\text{ with gradient } \leq\frac{t}{\bar{\rho}}\text{ and
}1\geq\inf_{x'\in\BB_{2\bar{\rho}}(x)}\frac{|x'-x|}{\text{dist}_{\SS}(x,x')}>\frac{9}{10}.$$
This means that, independently on $x$, $\S_{x,t}$ is a graph over
$T_x\SS$. Moreover, as shown in Fig. \ref{graph}, using the
Pythagorean theorem gives that
\begin{equation}\label{proj}
\text{the projection of } \S_{x,t}\text{ onto }T_x\SS \text{
contains }\DD_{\sqrt{t^2-\frac{t^4}{\bar{\rho}^2}}}(x).
\end{equation}
Furthermore, if $y\in B_t(x)\cap\Sigma$ and
$\text{dist}_\S(x,y)\geq 2t$ then $y$ cannot be in $\S_{x,t}$,
otherwise applying \eqref{eq2} gives
\begin{equation}\label{notins}
t\geq |y-x| > \frac{9}{10}\text{dist}_\S(x,y) \geq \frac{18}{10}t.
\end{equation}
$y$ is in a different component of $B_t(x)\cap\Sigma$. After
defining an orientation $y$ is either above or below $\S_{x,t}$.
For the same reason we can also add that
$\BB_{\bar{\rho}}(x)\cap\BB_{\bar{\rho}}(y)=\emptyset$.

Corollary \ref{criteria1} tells us that under certain conditions
regarding the orientation, if a CMC surface is a graph over
another CMC surface with the same constant mean curvature, then it
is almost stable. We are about to prove some lemmas which tell us
when that happens and how large the almost stable domain is.
This lemma shows how, if two pieces of $\S$ are close, then they
must be graphs over the same plane.
 \ble There exists
$\alpha_1>0$ so: For any $\alpha<\alpha_1$ and $x\in\Sigma'$ then
any component of $B_\alpha(x)\cap\SS'$ is a graph over $T_x\SS$.
\ele
\begin{proof}
Let us assume that there exists a component of
$B_\alpha(x)\cap\SS'$ which is not a graph over $T_x\SS$. Then
there exists $y\in B_\alpha(x)\cap\SS'$ such that $T_x\SS\perp
T_y\SS$. If $\alpha_1$ is small enough, it is clear from Fig.
\ref{graph2}
\begin{figure}[h]
    \setlength{\captionindent}{4pt}
    \begin{minipage}[t]{0.5\textwidth}
    \centering\input{fig6.pstex_t}
  \caption{}
   \label{graph2}
    \end{minipage}
\end{figure}
that
$\S_{x,\frac{\bar{\rho}}{2}}\cap\S_{y,\frac{\bar{\rho}}{2}}\neq\emptyset$
that is $y\in\BB_{\bar{\rho}}(x)$. How we have chosen $\bar{\rho}$
implies that $y$ must be part of a graph. Notice that we are also
using the fact that we are slightly away from the boundary.
$\S_{x,\frac{\bar{\rho}}{2}}\cap\S_{y,\frac{\bar{\rho}}{2}}$ could
be empty if one of the two sets reaches $\partial \S$ before they
intersect. How small $\alpha_1$ must be will depend also on
$\bar{\epsilon}$.

\end{proof}
In particular, it follows that if pieces of $\S'$ are very close
then not only are they graphs over the same plane, they are graphs
over each other. The idea is that if two graphs are almost flat
over two different planes but they cannot intersect, then if the
two graphs are close enough these two planes must have almost the
same slope. One of the two graphs can therefore be seen as a graph
over the other and this is what the next lemma is about.
\ble\label{harnack2} There exists $\alpha_2>0, C_1>0$ and $s>0$
so: Let $x,y\in\Sigma'$ such that $|x-y|\leq \alpha<\alpha_2$,
$d_{\Sigma}(x,y)\geq 2\alpha$ and $\lng n(x),n(y) \rng >0$ then
$\BB_s(y)$ contains a graph $\{z+u(z)n(z)\}$ over a domain
containing $\BB_{\frac{s}{4}}(x)$, $|\nabla u|+|u|\leq \alpha
C_1$. \ele
\begin{proof}
Assume $\alpha_2<\alpha_1$. We know from \eqref{notins} that $y$
is in a different component of $B_\alpha(x)\cap\Sigma'$ and that
$\BB_{\bar{\rho}}(x)\cap\BB_{\bar{\rho}}(y)=\emptyset$. If
$\alpha_2$ is sufficiently small, it now follows that
$\Sigma_{x,\bar{\rho}}$ and $\Sigma_{y,\bar{\rho}}$ contain two
graphs over the same plane, the smaller $\alpha_2$ is the bigger
the graph is. There exists $s>0$ such that $\BB_s(x)$ and
$\BB_s(y)$ contain respectively a graph $u_1$ and $u_2$ over
$\DD_\frac{s}{2}(x)$. $\lng n(x),n(y)\rng>0$ implies that the two
intrinsic disks have equal constant mean curvature and \eqref{CO}
gives $\sup_{\DD_{\frac{s}{4}(x)}}|u_1-u_2|\leq\alpha C_0$. The
function
\begin{equation}u(x)=\min\{t\in \mathbb{R}_+|x+tN(x)\in\BB_s(y)\}\end{equation}
is well defined over $\B_{\frac{s}{4}}(x)$ and $|\nabla u|+|u|\leq
\alpha C_1$.
\end{proof}
Lemma \ref{stability2} and  Lemma \ref{harnack2} give a better
criteria to find $\delta$-stable domain.
\begin{corr}\label{criteria2} Given $\delta>0$ there exists
$\alpha_3>0$ and $s>0$ so: Let $x,y\in\Sigma'$ such that
$|x-y|\leq \alpha<\alpha_3$, $d_{\Sigma}(x,y)\geq 2\alpha$ and
$\lng n(x),n(y)\rng>0$ then $\BB_{\frac{s}{4}}(x)$ is
$\delta$-stable.
\end{corr}
\begin{proof}
Let $\alpha_3<\min(\frac{\delta}{C_1},\alpha_2)$ and apply Lemma
\ref{stability2} and Lemma \ref{harnack2}.
\end{proof}
%
%
In sum, we have proven that when two points $x,y \in \Sigma'$ are
close enough to each other (Euclidean distance) and satisfy the
condition on the orientation $\langle N(x),N(y)\rangle>0$, then a
little neighborhood of each point is $\delta$-stable. We shall
notice that the closer two pieces are the smaller $\delta$ is. The
next step is to go from a little almost stable domain to a large
one.

If we need a very large $\delta$-stable geodesic ball, first of
all we need the geodesic ball to be contained in $\S$. In order to
achieve this we certainly cannot be anywhere in the surface but
sufficiently away from its boundary. If we move away from the
boundary, as long as the objects we are working with are contained
in $\S$, thanks to the Harnack inequality we can find conditions
that guarantee the existence of arbitrary large $\delta$-stable
domains. This is what we prove in the next lemma and corollaries.
In Fig. \ref{f:f2} it is shown how if two pieces of $\S^1$ are
close then their extensions will have to stay relatively close.
\begin{figure}[htbp]
    \setlength{\captionindent}{4pt}
    \begin{minipage}[t]{0.5\textwidth}
    \centering\input{fig7.pstex_t}
   \caption{}
   \label{f:f2}
    \end{minipage}
\end{figure}
\begin{lemma}\label{largestab1}
For each $0<l<\bar{l}-(1+\bar{\epsilon})$ there exist $\alpha_l>0$
and $C_l>0$ so: Given $\alpha<\alpha_l$ and $x,y\in\Sigma^1$ such
that $|x-y|\leq \alpha$, $d_{\Sigma}(x,y)\geq C_l $ and $\lng
n(x),n(y) \rng >0$ then for each $x'\in \BB_l(x)$ there exists
$y'\in \S'$ such that $|x'-y'|\leq \alpha_1$,
$d_{\Sigma}(x',y')\geq 2\alpha_1$ and $\lng n(x'),n(y') \rng >0$.
\end{lemma}
\begin{proof}
Fix $N\in\mathbb{N}$ such that $\frac{4l}{s}\leq
N<\frac{4l}{s}+1$, $s$ as in Corollary \ref{criteria2} and assume
$\alpha_l<\min(\alpha_2,\frac{\alpha_2}{C_0^N})$, $C_0$ being the
constant as in \eqref{CO}. Our goal is to find $C_l$, note that
$\BB_l(x)\subset\S'$. Let $x'\in\BB_l(x)$ then there exists a
geodesic $\gamma(t),\SP t\in [0,1]$ such that $\gamma(0)=x$,
$\gamma(1)=x'$ and $\text{length}(\gamma)\leq l$. Fix a partition
$\mathcal{Q}$ of $[0,1]$, $\mathcal{Q}=\{t_i\in [0,1]| 0\leq i\leq
T \}$, such that
\begin{equation}
\begin{cases}
t_0=0, t_T=1\\
d_\Sigma(\gamma(t_i)=x_i,\gamma(t_{i+1})=x_{i+1})\leq
\frac{s}{4}\\
T \leq N
\end{cases}.\end{equation}
Since $|x-y|\leq \alpha\leq\alpha_2$ and $d_{\Sigma}(x,y)\geq
C_l\geq 2\alpha_2$ then Lemma \ref{harnack2} gives $\BB_s(x)$ and
$\BB_s(y)$ contain respectively a graph $u_1$ and $u_2$ over
$\DD_{\frac{s}{2}(x)}$ and
$\sup_{\DD_{\frac{s}{4}(x)}}|u_1-u_2|\leq C_0 \alpha$. Therefore,
let $z_1\in \DD_{\frac{s}{4}}(x)$ such that $x_1=u_1(z_1)$ and let
$y_1=u_2(z_1)$, then $|x_1-y_1|\leq C_0 \alpha\leq\alpha_2$,
$d_{\Sigma}(x_1,y_1)\geq
d_{\Sigma}(x,y)-d_{\Sigma}(x_1,x)-d_{\Sigma}(y_1,y)\geq
C_l-\frac{s}{4}-\frac{5s}{4}\geq 2\alpha_2$ if $C_l$ is big
enough, and $\lng N(x_1),N(y_1)\rng>0$. As long as
$d_{\Sigma}(x_i,y_i)\geq 2\alpha_2$ we can apply Lemma
\ref{harnack2}. We can repeat this argument $N$ times as long as
$C_l-N\frac{3s}{2}\geq 2\alpha_2$.
\end{proof}
\begin{corr}\label{largestab3}
For each $0<l<\bar{l}-(1+\bar{\epsilon})$ there exist $\alpha_l>0$
and $C_l>0$ so: Given $\alpha<\alpha_l$ and $x,y\in\Sigma^1$ such
that $|x-y|\leq \alpha$, $d_{\Sigma}(x,y)\geq C_l $ and
$<n(x),n(y)>>0$ then for each $x'\in \BB_l(x)$, $\BB_t(y)$
contains a graph $\{z+u(z)n(z)\}$ over a domain containing
$\BB_{\frac{t}{4}}(x)$, $|\nabla u|+|u|\leq \alpha C_1$.
\end{corr}
\begin{proof}
Apply Lemma \ref{largestab1} and Lemma \ref{harnack2}.
\end{proof}
\begin{corr}\label{largestab2}
For each $0<l<\bar{l}-(1+\bar{epsilon})$ and $\delta>0$ there
exist $\alpha_{l,\delta}>0$ and $C_l>0$ so: Given
$\alpha<\alpha_{l,\delta}$ and $x,y\in\Sigma^1$ such that
$|x-y|\leq \alpha$, $d_{\Sigma}(x,y)\geq C_l$ and $\lng
N(x),N(y)\rng>0$ then $\BB_l(x)$ is $\delta$-stable.
\end{corr}
\begin{proof}
Let $\alpha_{l,\delta}<\min(\alpha_l,\frac{\delta}{C_1})$,
according to Corollary \ref{criteria1}, we need to find a $u$
which is a positive solution of the CMC graph equation over
$\BB_l(x)$. Corollary \ref{largestab3} gives that the latter is
true locally. In fact, fix $x_i\in\BB_l(x)$ such that $\BB_s(x_i)$
and $\BB_{\frac{s}{4}}(x_i)$ are both finite coverings for
$\BB_l(x)$. From Lemma \ref{largestab1} it follows that
$\BB_s(x_i)$ contains a graph $\{z+u_i(z)n(z)\}$ over a domain
containing $\BB_{\frac{s}{4}}(x_i)$, $|\nabla u_i|+|u_i|\leq
\alpha C_1$. The function $u(y):=u_i(y)$ if
$y\in\BB_{\frac{s}{4}}(x_i)$ is a well defined function over
$\BB_l(x)$ such that $|\nabla u|+|u|\leq \alpha C_1<\delta$.
Applying Corollary \ref{criteria1} gives this corollary.
\end{proof}

%
%
%
%

\newsection{The non-existence of large almost-stable
domains and the uniform bound}

We have seen when it happens that $\S$ contains a large almost
stable domain. In this section we show that an almost stable
domain cannot be too large. Using these two facts together we
prove a uniform bound on the number of graphs.

In order to continue with this proof by contradiction we state the
following result by Sirong Zhang \cite[\bf Theorem 0.1.]{Zh}, and
the Bishop Volume Comparison Theorem \cite[\textbf{Theorem
1.3.}]{ScY}:
\begin{theorem}\label{Zh} There exists a $C$ such that given any $l>0$
there exists an $h>0$ so: If $\BB_l(0)$ is a "constant mean
curvature equal to $h$", $\delta$-stable intrinsic disk with
trivial normal bundle then $\sup_{\BB_{\frac{l}{2}}(0)}|A|^2\leq
\frac{C}{l^2}$.
\end{theorem}

\begin{theorem}[Bishop Volume Comparison Theorem]\label{bishop}
Let $M$ be an $n$-dimensional complete Riemannian manifold with
$Ric(M)\geq(n-1)K$. Then for any $x\in M$ and $R>0$,
$\frac{Vol(\B_R(x))}{V(K,R)}$ is a non-increasing function in $R$.
Hence, $$Vol(\B_R(x))\leq V(K,R),$$ where $V(K,R)$ is the volume
of the geodesic ball of radius $R$ in the space form $M_K$.
\end{theorem}
Theorem \ref{Zh} can be thought as a generalization of
\cite{doCPe,FCSch} and we will see how it is essentially what
determines how big $\bar{l}$ is. Our surface has trivial normal
bundle since it is orientable and this will be proved later in
Proposition \ref{finale}. %
Theorem \ref{bishop} gives that for any $x\in\S$
\begin{equation}
Vol(\B_R(x))\leq V(G,R), \quad\quad G \text{ as in \ref{bishop1}}.
\end{equation}
The following proposition uses what we have proved in the previous
sections and Theorem \ref{Zh} to show that $\S'$ does not contain
a large almost stable domain. Roughly speaking if we take $l$
large and assume that $\BB_l(x)\subset\Sigma^1$ is
$\delta$-stable, Theorem \ref{Zh} implies that
$\BB_\frac{l}{2}(x)$ is almost flat. This forces the intrinsic
disk to leave the unit ball.
\begin{proposition}\label{nonexist} Given $\delta>0$ there exists
$l_{\delta}>0$ so: For any $l\geq l_{\delta}$ and $x\in \Sigma^1$
if $\BB_l(x)$ is $\delta$-stable then it is not contained in $
\Sigma^1$.
\end{proposition}
\begin{proof}
Let us fix $l_{\delta}>\max(\frac{20}{9},\sqrt{\frac{C}{\rho}})$.
Let $C$ be as in Theorem \ref{Zh} and $\rho$ as in \eqref{eq2}.
Being $\BB_l(x)$ $\delta$-stable Theorem \ref{Zh} implies
$$\sup_{\BB_{\frac{l}{2}}(x)}|A|^2\leq \frac{C}{l^2} \leq \rho.$$
\eqref{eq2} implies that
$$\inf_{\BB_{\frac{l}{2}}(x)}\frac{|x-y|}{d_{\Sigma}(x,y)}\geq
\frac{9}{10}.$$ Taking $y\in\BB_{\frac{l}{2}}(x)$ such that
$d_{\Sigma}(x,y)>\frac{10}{9}$, we have $|x-y|> 1$. This proves
that $\BB_l(x)$ is not contained in $\S^1$.
\end{proof}
We have proved so far that decreasing the Euclidean distance
between two points gives a large $\delta$-stable domain as long as
we increase their intrinsic distance. In the following lemma we
apply the Bishop Volume Comparison Theorem and a lower bound on
the area of each piece to prove that the more graphs there are in
a small ball the larger the intrinsic distance
becomes. This is what Fig. \ref{f:f3} illustrates.%
\begin{figure}[htbp]
    \setlength{\captionindent}{4pt}
    \begin{minipage}[t]{0.5\textwidth}
    \centering\input{fig8.pstex_t}
  \caption{}
   \label{f:f3}
    \end{minipage}
\end{figure}

Let us fix $\delta, \bar{\epsilon}>0$ small, let $l_1=l_{\delta}$
as given by Proposition \ref{nonexist} and let $\bar{l}$ in
Theorem \ref{main2} to be equal to $l_1+1+\bar{\epsilon}$. In
other words,
\begin{quote}
$\S^1$ does not contain a $\delta$-stable geodesic ball of radius
bigger than $l_1$.
\end{quote}
Let us fix $0<r<\alpha_{l_1,\delta}$ where $\alpha_{l_1,\delta}$
is taken as in Corollary \ref{largestab2}. This means that
\begin{equation}\label{radius1}
\text{For any }x\in\S'\text{, }B_r(x)\cap\S'\text{ consists of
graphs over }T_x\S
\end{equation}
and also that
\begin{equation}\label{radius2}
\begin{array}{l}
\text{There exists }C_{l_1}>0\text{ such that if }x,y\in
B_r(x)\cap\S^1\text{, }d_{\Sigma}(x,y)\geq C_{l_1}\\
\text{ and }\lng N(x),N(y)\rng>0\text{ then }\BB_{l_1}(x)\text{ is
}\delta\text{-stable}.\end{array}
\end{equation}

Given $x\in\S^1$ let $n_x$ be the number of components of
$B_r(x)\cap\S^1$. The area of each component $\S_i^x$ could go to
zero if they accumulate toward the boundary of the ball.
Nonetheless we have proven that these graphical pieces continue
outside the ball, \ref{proj}. Therefore we have a uniform lower
bound on the area of $\S_i^x$; in other words
\begin{equation}
\text{There exits }\epsilon>0 \text{ such that
}\text{Area}(\S_i^x)>\epsilon \text{ independently of }x\text{ and
 }i.
\end{equation}
\ble\label{lambda} Given $\lambda>0$ there exists $n_\lambda>0$
so: If $x\in\S^1$ and $n_x>n_\lambda$ then there exist
$y,y'\in\B_r\cap\S^1$ such that $\text{dist}_\S(y,y')>\lambda$ and
$\lng n(y),n(y')\rng>0$.\ele
\begin{proof}
Theorem \ref{bishop}, that is, Bishop Volume Comparison Theorem,
gives an upper bound for the area of $\BB_\lambda(x)$, namely
$$\text{Area}(\BB_\lambda(x))<V(G,\lambda).$$ At this point it follows easily that
if $n_x\in\mathbb{N}$, $n_x>\frac{V(G,\lambda)}{\epsilon}$ then
there exists $y_1\in\B_r\cap\S^1$ which is not in
$\BB_\lambda(x)$, i.e. $\text{dist}_\S(x,y_1)>\lambda$. $\lng
n(x),n(y_1)\rng>0$ does not necessarily happen. Take $\bar{V}\in
\mathbb{N}$ such that
$\bar{V}-1<\frac{V(G,\lambda)}{\epsilon}\leq\bar{V}$ and let
$n_\lambda=\bar{V}^{\bar{V}}$. If $n_x\geq n_\lambda$ then there
exist at least $\bar{V}$ distinct $y_i,\SP i:=1,...\bar{K}$, in
different component of $B_r(x)\cap\S^1$ such that
$\text{dist}_\S(x,y_i)>\lambda$. Fixed $y_1$ there exists $y_2$
among the $y_i$ such that $\text{dist}_\S(y_1,y_2)>\lambda$. At
this point either $\lng n(x),n(y_1) \rng >0$ or $\lng n(x),n(y_2)
\rng >0$ or $\lng n(y_1),n(y_2)\rng>0$.
\end{proof}
The following corollary uses Proposition \ref{nonexist} and Lemma
\ref{lambda} to obtain the upper bound on the number of graphs.
\begin{corr}
For any $x\in\S^1$, $n_x\leq n_{C_{l_1}}$.
\end{corr}
\begin{proof}
If $n_x> n_{C_{l_1}}$ then Lemma \ref{lambda} with $\lambda$ equal
to $C_{l_1}$ gives that there exist $y,y'\in\B_r(x)\cap\S^1$ such
that $\text{dist}_\S(y,y')>C_{l_1}$ and $\lng n(y),n(y')\rng>0$.
Using Lemma \ref{largestab2} gives that $\BB_{l_1}(x)$ is
$\delta$-stable and, by Proposition \ref{nonexist}, cannot be
contained in $\S^1$. Since $\text{dist}_\S(y,y')>C_{l_1}>l_1$,
$y'$ is not in $\S^1$. This gives the contradiction that implies
$n_x\leq n_{C_{l_1}}$.
\end{proof}

\newsection{Multi-valued graphs in CMC surfaces}

In this final section we show that the $\S^1_n$ converges $C^2$ to
an embedded minimal disk $\S_\infty$ that contains a multi-valued
graph. It follows that the CMC surfaces in the sequence contain a
multi-valued graph as well. The limit surface is embedded and
minimal by a standard argument which will be sketched below. To
prove that it is simply connected we need more work and well-known
topological results.

Let $r>0$ be as defined in \eqref{radius1} and \eqref{radius2},
that is, there exists a finite covering $B_r(x^n_i)$ for $\S^1_n$
where everything is graphical over $T_{x^n_i}\S_n$ and the number
of graphs is uniformly bounded. We can also assume that the number
of balls involved is uniformly bounded with respect to $n$. Going
to a subsequence, we can assume $x^n_i$ converging to $x_i$ and
$T_{x^n_i}\S_n$ converging to a certain $T_{x_i}\S_\infty$. Using
the argument outlined in Section \ref{proof}, the fact that the
number of graphs is uniformly bounded and the maximum principle
for minimal surfaces gives that the limit is an embedded minimal
surface.
\begin{figure}[h]
    \setlength{\captionindent}{4pt}
    \begin{minipage}[t]{0.5\textwidth}
    \centering\input{fig9.pstex_t}
  \caption{}
  \label{emb}
    \end{minipage}
\end{figure}
Fig. \ref{emb} illustrates the two types of intersection that
could occur if the limit is not embedded: A cross intersection,
type A, and a tangential intersection, type B. However type A
cannot be a continuous limit of embedded surfaces. Type B, which
could be the limit of a sequence of embedded surfaces, cannot
occur because of the maximum principle for minimal surfaces. By
continuity the curvature of this minimal surface is large at zero.

To prove that $\S_\infty$ is simply connected we use some results
about Jordan curves \cite{GuPo} and the following theorems:

\begin{theorem}\label{orient}\cite[\textbf{Chapter 3}]{GuPo}
Every compact hypersurface in Euclidean space without boundary is
orientable.
\end{theorem}
\begin{theorem}\label{orient2}\cite[\textbf{Chapter 3}]{GuPo}
If $M$ is an orientable surface, so is $M$ minus one point.
\end{theorem}
\begin{theorem}\label{tor}\cite[\textbf{Corollary 3.28.}]{Ha}
If $M$ is a closed connected $n$-manifold, the torsion subgroup of
$H_{n-1}(M,\mathbb{Z})$ is trivial if $M$ is orientable and
$\mathbb{Z}_2$ if $M$ is not orientable.
\end{theorem}
In our case $\S_\infty$ is a closed connected 2-manifold, and it
is a consequence of Theorem \ref{tor} that, if it is orientable,
the torsion subgroup of its fundamental group is trivial. The next
proposition shows that $\S_\infty$ is the embedded minimal disk we
have been looking for.
\bprop \label{finale} $\S_\infty$ is an embedded simply connected
minimal surface such that, $0\in\Sigma_\infty\subset
B_{1}\subset\R^3$, $\partial \Sigma_\infty \subset\partial
B_{1}$, $\sup_{\Sigma\cap B_{1}}|A|^2\leq 4C^2$ and
$|A|^2(0)=C^2$. \eprop
\begin{proof}
The conditions on the second fundamental form, that is
$$\sup_{\Sigma\cap B_{1}}|A|^2\leq 4C^2 \text{ and } |A|^2(0)=C^2,$$ are
a consequence of the $C^2$ convergence. We have already proved
that it is embedded and what we are left to prove is that
$\S_\infty$ is simply connected. Let us prove that it is
orientable first.
\begin{figure}[htbp]
    \setlength{\captionindent}{4pt}
    \begin{minipage}[t]{0.5\textwidth}
    \centering\input{fig10.pstex_t}
   \caption{}
\label{boundary}
    \end{minipage}
\end{figure}
We want to prove that $\S_\infty$ is
homeomorphic to a compact embedded surface minus a finite number
of points. Because it is an embedded minimal surface, $\partial
\S_\infty$ is a finite number of disjoint loops $\xi_i$,
$i:=1,...,I$ which do not have self intersections, Fig.
\ref{boundary}. These loops lie on $\partial B_{1}(0)$ minus one
point, hence we have essentially a finite number of Jordan curves
in the plane. We can glue $I$ disks to $\S_\infty$ in a way that
the result is an embedded compact surface: Each Jordan curve
$\xi_i$ divides the plane into an inner and an outer region, and
can be thought of as the boundary of a simply connected domain,
namely a disk $D_i$. If a loop $\xi_i$ lies in the inside of
another loop $\xi_j$ then we lift $D_j$ so that it does not
intersect $D_i$. Fig. \ref{patch} illustrates how we are gluing
these disks to the surface.
\begin{figure}[htbp]
    \setlength{\captionindent}{4pt}
    \begin{minipage}[t]{0.5\textwidth}
    \centering\input{fig11.pstex_t}
   \caption{}
   \label{patch}
    \end{minipage}
\end{figure}
Since the number of loops is finite we repeat this a finite number
of times and obtain in the end a new surface
$$\bar{\S}=\{\S_\infty \bigcup_{i=1}^ID_i\} / \sim$$ where $\sim$
is the relation that identify $\xi_i$ with $\partial D_i$.
$\bar{\S}$ is a compact embedded surface without boundary and
therefore orientable by Theorem \ref{orient}. Theorem
\ref{orient2} implies that $\bar{\S}$ take out a finite number of
points is still orientable, that is $\bar{\S}\backslash
\bigcup_{i=1}^ID_i=\S_\infty$. Theorem \ref{tor} tells that
$\pi_1(\S_\infty)$ is torsion free.

Let $\gamma :\SP S^1\longrightarrow \S_\infty$ be a closed path
and
 $B_\sigma(\gamma(t_i))$ a finite covering for $\gamma$ such that $\sigma<\min(1-|\gamma|,r)$, and
the number of components of $B_\sigma(\gamma(t))$ is
non-increasing. This is possible after going to a subsequence,
assuming $n$ large because of the uniform bound. Fix a starting
point $\gamma(t_0)$, an orientation on $\gamma$, and let
$$\gamma_n(t_0):=\gamma(t_0)+s_0N_{\gamma(t_0)}\in\S^1_n.$$ Moving continuously on
$\gamma(t)$ we obtain a new path
$$\gamma_n(t):=\gamma(t)+s(t)N_{\gamma(t)}\in\S^1_n.$$
The conditions on $\sigma$ and $n$ force the path to close up
after it moves around $\gamma$ a finite number of times,
$k\in\mathbb{N}$. Since $\S^1_n$ simply connected, there exists a
map $\Gamma_n:\SP D_1\longrightarrow \S^1_n$ such that
$\Gamma_n|_{S^1}=\gamma_n$. Define $\Gamma:\SP D_1\longrightarrow
\S_\infty$,
$$\Gamma(x):=\lim_{n\longrightarrow\infty}\Gamma_n(x).$$
The existence of this map $\Gamma$ proves that $k\gamma$ is
homotopic to a point. Since $\pi_1(\S_\infty)$ is torsion free,
this implies that $\gamma$ itself is homotopic to a point. Since
$\gamma$ could be any path on $\S_\infty$ we have proved that
$\S_\infty$ is simply connected.
\end{proof}
Finally we prove that $\S^1_n$ and therefore $\S_n$ contains a
multi-valued graph. In Proposition \ref{finale} we proved that
$\S_\infty$ is an embedded simply connected minimal surface such
that, $$0\in\Sigma_\infty\subset B_{1}\subset\R^3, \partial
\Sigma_\infty \subset \partial B_{1}, \sup_{\Sigma\cap
B_{1}}|A|^2\leq 4C^2=4|A|^2(0).$$ Taking $C=C(N,\omega,\epsilon)$
as in Theorem \ref{result}, the same theorem gives that
$\S_\infty$ contains and $N$-valued graph. Let $u$ be this
$N$-valued graph, defined over $\{(\rho ,\theta )| r\leq \rho \leq
s,|\theta |\leq N\pi \}$ as described in Definition
\ref{multigraph}. This is how we build an $N$-valued graph in
$\S_n$: Given $r\leq \bar{\rho} \leq s$, define
$u_{\bar{\rho}}(\theta)=u(\bar{\rho},\theta)$. Consider
$u_{\bar{\rho}}$ as a path on $\S_\infty$ starting at
$u_{\bar{\rho}}(-N\pi)$. Assuming $n$ large, $\S^1_n$ moves closer
and closer to $\S_\infty$ and there exists a continuous function
$\phi(\theta)$ such that
$$u^n_{\bar{\rho}}(\theta)=u_{\bar{\rho}}(\theta)+\phi(\theta)\textbf{e}_3\in\S^1_n$$
is well defined. The function $u^n(\rho,\theta)=u^n_\rho(\theta)$
defined over $\{(\rho ,\theta )| r\leq \rho \leq s,|\theta |\leq
N\pi \}$ is an $N$-valued graph.

Notice that as $\S^1_n$ moves closer and closer to $\S_\infty$,
$\S^1_n$ and $\S_\infty$ are "parallel surfaces." Not only does
$\S_n$ contain an $N$-valued graph, but the properties of this
graph, such as the upper bound on the gradient, are preserved.

 \renewcommand{\theequation}{A-\arabic{equation}}
 \renewcommand{\thesubsection}{A-\arabic{subsection}}
 \renewcommand{\thetheorem}{A-\arabic{theorem}}
 \setcounter{equation}{0}  
 \setcounter{theorem}{0}

 \section*{APPENDIX A}

In this appendix we provide examples of CMC surfaces containing
arbitrary large multi-valued graphs. We use the method of
successive approximations to build a sequence of normal variations
of the helicoid that converges to an embedded and simply connected
CMC surface containing a multi-valued graph.

Let $\S_h=\{x+u(x)N_{\S}(x),x\in\S\}$ be a normal variation of
$\S$, where $\S$ is any minimal surface. $\S_h$ is a CMC surface
with mean curvature equal to $H$ if $u(x)$ satisfies the following
equation \cite{Ka1, Ka2,Nit}:
\begin{equation}\label{solution}
Lu=H+Q(u), \quad \quad \text{where } Lu=\Delta u +|A|^2u
\end{equation}
is the linearized operator. $Q$ is a quadratic and higher order
function in $u,\SP u_i,\SP u_{ij}$ where $i,\SP j\in\{1,2\}$, with
geometric invariants of $\Sigma$ as coefficients. Before we prove
the existence of a constant mean curvature normal variation of the
helicoid we need to describe some properties of the function $Q$.

Let $C^{k,\lambda}(\Sigma)$ be the standard subset of
$C^k(\Sigma)$ consisting of functions whose $k$-th partial
derivatives are H\"{o}lder continuous with exponent $\lambda$ in
$\Sigma$ and let $\|\cdot\|_{k,\lambda}$ be the notation for the
H\"{o}lder norm. Let us define $\Omega_\delta$, subset of
$C^{2,\lambda}(\Sigma)$, in the following way:
\begin{equation} \Omega_\delta:=\{u\in C^{2,\lambda}:
\|u\|_{2,\lambda}<\delta\}.
\end{equation}

The following lemma follows from \cite[Lemma C.2]{Ka1} and it is a
consequence of the fact that $Q$ is a quadratic and higher order
function. It says that $\|Q(u)\|_{0,\lambda}$ decays faster than
$\|u\|_2\|u\|_{2,\lambda}$.
\begin{lemma}\label{ibound}
There exist $\delta_1>0$ and $C_1>0$ so: If $u \in
C^{2,\lambda}(\Sigma)$ and $|A||u|,|u_i|<\delta_1$ then
$$\|Q(u)\|_{0,\lambda}<C_1\|u\|_2\|u\|_{2,\lambda}.$$
\end{lemma}
As a consequence of Lemma \ref{ibound}, we have a new corollary
that relates $\|Q(u)\|_{0,\lambda}$ and $\|u\|_{2,\lambda}$:
\begin{corr}\label{iibound}
Given $1>C_2>0$ there exists $\delta_2>0$ so: If $u\in
\Omega_{\delta_2}$ and $|A||u|<\delta_2$ then
$$\|Q(u)\|_{0,\lambda}<C_2\|u\|_{2,\lambda}.$$
\end{corr}
\begin{proof}
Let $\delta_2<\min(\frac{C_2}{C_1},\delta_1)$. This implies that
$C_1\|u\|_2<C_1\|u\|_{2,\lambda}<C_1\delta_2<C_2$ and that
$|A||u|,|u_i|<\delta_2<\delta_1$. Therefore, we can apply
\ref{ibound} and we have
$$\|Q(u)\|_{0,\lambda}<C_1\|u\|_2\|u\|_{2,\lambda}<C_2\|u\|_{2,\lambda}.$$
\end{proof}

Let $\Sigma$ be a simply connected disk in the helicoid that
contains a multi-valued graph. Due to the domain monotonicity and
continuity of eigenvalues \cite{Ch} we can also assume that $0$ is
not an eigenvalue for $L$ on $\Sigma$ and therefore that the
Dirichlet
\begin{equation}\label{luw}
\begin{cases}
Lu=w \\ u |_{\partial\S}=0
\end{cases}
\end{equation}
has a unique solution $u\in C^{2,\lambda}(\S)$ for any $w\in
C^{0,\lambda}(\S)$. Assuming that $0$ is not an eigenvalue for $L$
gives also the following lemma \cite[Theorem 5.3 and page
109]{GiTru}:
\begin{lemma}\label{estim}
There exists a constant $B$ depending only on $\S$ so: Let $u\in
C^{2,\lambda}(\S)$ be the unique solution for \eqref{luw} then
\begin{equation}
\|u\|_{2,\lambda}<B\|w\|_{0,\lambda}.
\end{equation}
\end{lemma}

We will prove that there exists $H>0$ such that a solution $u$ for
the Dirichlet  problem
\begin{equation}
\begin{cases}
Lu=H+Q(u)\\
u|_{\partial\S}=0
\end{cases}
\end{equation}
exists and $\|u\|_{L^\infty}$ is small. The existence of a fixed
neighborhood of the helicoid where the normal exponential map is
injective guarantees that $\S_H=\{x+u(x)N_{\S}(x),x\in\S\}$ is
embedded, if $\|u\|_{L^\infty}$ is small enough. What we are about
to show is that if $H$ is small enough we can build a sequence of
normal variations $u^n$ of the helicoid that converges to a CMC
normal variation $u$. We will also show that $\|u\|_{L^\infty}$
can be as small as we want and consequently that the CMC normal
variation is embedded.

Let $u^1$ be the unique solution for
\begin{equation}
\begin{cases}
Lu^1=H \\
u^1 |_{\partial\S}=0
\end{cases}\end{equation}
and $u^n$ be the unique solution for
\begin{equation}
\begin{cases}
Lu^n=H+Q(u^{n-1}) \\
 u^n |_{\partial\S}=0
\end{cases}.\end{equation}
Lemma \ref{estim} implies that
$$\|u^1\|_{2,\lambda}<BH$$
and also that
\begin{equation}\label{estim2}
\|u^k\|_{2,\lambda}<B(H+\|Q(u^{k-1})\|_{0,\lambda}).
\end{equation}
The existence of a solution for
$$Lu=H+Q(u)$$
will follow clearly, and we will see how, applying Arzela-Ascoli
to the sequence $u^n$ if we prove that there exists a constant $K$
such that $\|u^n\|_{2,\lambda}<K$ uniformly in $n$. Fix $C_2$ in
Corollary \ref{iibound} so that $\epsilon=C_2B<1$, $B$ as in Lemma
\ref{estim}. We will prove by strong induction that
\begin{quote}
if $H$ is so that $BH(1+\frac{1}{1-\epsilon})<\delta_2$ then
$u^n\in\Omega_{\delta_2}$ for any $n$,\end{quote} that is what we
wanted.
We have already that
$$\|u^1\|_{2,\lambda}<BH<\delta_2,$$ namely the statement is true for $n=1$.
Let us prove that
\begin{quote}
"true for $n=1$ implies true for $n=2$." \end{quote}
Assuming $\|u^1\|_{2,\lambda}<\delta_2$, we can apply Lemma
\ref{ibound} that gives that
$\|Q(u^1)\|_{0,\lambda}<C_2\|u^1\|_{2,\lambda}<\epsilon H$,
therefore
\begin{equation}
\|u^2\|_{2,\lambda}\le B(H+\|Q(u^1)\|_{0,\lambda}\le B(H+\epsilon
H)<\delta_2.
\end{equation}
Let us prove that
\begin{quote}
"true for all $k$ with $k\leq n$ implies true for $k=n+1$."
\end{quote}
"True for all $k$ with $k\leq n$" means that
$\|u^k\|_{2,\lambda}<\delta_2$ for $k\leq n$ and therefore Lemma
\ref{ibound} gives that
\begin{equation}\label{induction}
\|Q(u^{k})\|_{0,\lambda}<C_2\|u^{k}\|_{2,\lambda}<C_2B(H+\|Q(u^{k-1})\|_{0,\lambda})<
\epsilon (H+\|Q(u^{k-1})\|_{0,\lambda})
\end{equation}
for $k\leq n$. Applying \eqref{induction} $n$ times we have
\begin{multline}
\|u^{n+1}\|_{2,\lambda}< B(H+\|Q(u^n)\|_{0,\lambda})< B(H+\epsilon
(H+\|Q(u^{n-1})\|_{0,\lambda}))< \\<B(H+\epsilon (H+\epsilon
(H+\|Q(u^{n-2})\|_{0,\lambda})))<\\<BH(1+\sum_{k=1}^n\epsilon^k)<BH(1+\frac{1}{1-\epsilon})<\delta_2.
\end{multline}

Now that we have proved that $\|u^n\|_{2,\lambda}<\delta_2$
uniformly in $n$, using Arzela-Ascoli we can extract a subsequence
that converges $C^2$ to a certain $u\in C^2(\Sigma)$. Taking the
limit as $n$ goes to infinity on both sides of the equation
$$Lu^n=H+Q(u^{n-1})$$ gives that $$Lu=H+Q(u).$$ $u$ is therefore a
constant mean curvature normal variation of the helicoid. It is
clear from the proof that taking $H$ small gives
$\|u\|_{L^\infty}$ small. Consequently, the constant mean
curvature normal variation that we have built is also embedded.

\renewcommand{\theequation}{B-\arabic{equation}}
 \renewcommand{\thesubsection}{B-\arabic{subsection}}
 \renewcommand{\thetheorem}{B-\arabic{theorem}}
 \setcounter{equation}{0}  
 \setcounter{theorem}{0}
\section*{APPENDIX B}  

In this appendix we want to show that Theorem \ref{CoMi} follows
from Theorem \ref{main2} by rescaling. The result is true even
when the mean curvature is large but on a smaller ball. In other
words, surfaces with large constant mean curvature have tiny
multi-valued graphs around the origin. We give the idea of how
that happens when dealing with simple graphs.

Let us assume that the CMC surface $\Sigma$ is given as a graph
$u$ over the unit disk (given that it contains a multi-valued
graph this cannot actually happen globally but it is always
possible locally and we are assuming it now just for the sake of
simplicity). Then the new surface $\Sigma'$ given by
$w(x)=Ru(\frac{x}{R})$ defined over $D_R(0)$ is still a CMC
surface. In fact
 \beq
H(\Sigma')=\text{div}\left (\frac{\nabla w}{\sqrt{1+|\nabla
w|^2}}\right)=\frac{1}{R}H(\Sigma). \eeq
Therefore, assuming for instance $R<1$ the new CMC surface has
bigger mean curvature and the multi-graph happens in a smaller
ball (very tiny if $H$ is big). Rescaling preserves the existence
of the multi-valued graph but changes the hypotheses regarding
$|A|^2$. Since $k'_i=\frac{1}{R}k_i$, we have
$|A'|^2=\frac{|A|^2}{R^2}$ and therefore
$$|A'(0)|^2=\frac{|A(0)|^2}{R^2}=\frac{C^2}{R^2}\text{ and } \sup_{\Sigma'\cap B_R}|A'|^2=\frac{\sup_{\Sigma\cap
B_1}|A|^2}{R^2}\leq\frac{4}{R^2}C^2$$ hence Theorem \ref{CoMi}.

\end{document}